\documentclass[12pt,reqno]{amsart}
\usepackage{amsmath, amssymb, latexsym, amscd, amsthm,amsfonts,amstext}
\usepackage[mathscr]{eucal}
\usepackage{graphics}
\newcommand{\R}{\mathbb{R}}
\newcommand{\C}{\mathbb{C}}
\newcommand{\D}{\mathbb{D}}

\newcommand{\Z}{\mathbb{Z}}

\theoremstyle{plain}
\newtheorem{lemma}{Lemma}[section]
\newtheorem{theorem}[lemma]{Theorem}
\newtheorem{corollary}[lemma]{Corollary}

\newtheorem{definition}[lemma]{Definition}

\newtheorem{example}[lemma]{Example}

\newenvironment{proof*}{\vskip 2mm\noindent {}}{\hfill $\Box$ \vskip 2mm}

\begin{document}

\title[On hyperbolicity and tautness modulo an analytic subset...]{On hyperbolicity and tautness modulo an analytic subset of Hartogs domains}
\thanks{The research of the authors is supported by an NAFOSTED grant of Vietnam.}
\author[D. D. Thai, P. J. Thomas, N. V. Trao, M. A. Duc]{Do Duc Thai*, Pascal J. Thomas**, Nguyen Van Trao* \& Mai Anh Duc*}

\address{*Department of Mathematics \newline
Hanoi National University of Education\newline
136 XuanThuy str., Hanoi, Vietnam}

\email{ducthai.do@gmail.com}
\email{ngvtrao@yahoo.com}
\email{ducphuongma@gmail.com}

\address{**Universit\'e de Toulouse\\ UPS, INSA, UT1, UTM \\
Institut de Math\'e\-matiques de Toulouse\\
F-31062 Toulouse, France}

\email{pascal.thomas@math.univ-toulouse.fr}

\subjclass[2000]{Primary 32M05; Secondary 32H02, 32H15, 32H50.}
\keywords{Hyperbolicity modulo an analytic subset, tautness modulo an analytic subset, Hartogs domains.}

\begin{abstract} 
Let $X$ be a complex space and $H$ a positive homogeneous plurisubharmonic function $H$ on $X\times\C^m$.
Consider the Hartogs-type domain  $\Omega_{H}(X):=\{(z,w)\in X\times \C^m:H(z,w)<1 \}$. Let $S$ be an analytic subset of $X$.
We give
necessary and sufficient conditions for hyperbolicity and  tautness modulo
$S\times \C^m$ of $\Omega_{H}(X)$, with the obvious corollaries for the special case of
Hartogs domains.
\end{abstract}

\maketitle  
 
\section{Introduction}

Let $X$ be a complex space and $\varphi: X\to [-\infty,\infty)$ be an upper-semicontinuous function on $X$. The Hartogs domain
 $$\Omega_{\varphi}(X):=\{(x,z)\in X\times \C:|z|<e^{-\varphi(x)}\}$$
is a classical object in Several Complex Variables. In particular, in the past ten years, much attention has been given to the properties of Hartogs domains from the viewpoint of hyperbolic complex analysis. For instance, in  \cite{TD}, the authors obtained necessary and sufficient conditions for the hyperbolicity and the tautness of $\Omega_{\varphi}(X)$.  We refer readers to the articles \cite{TT},  \cite{DT}, and references therein for the development of related subjects. 

More generally if $H: X\times \C^m \to [-\infty,\infty)$
is an upper semicontinuous function such that $H(z,w) \ge 0, H(z,\lambda w)=|\lambda|H(z,w), \lambda \in \C, z \in X, w \in \C^m$, we put
$$
\Omega_H(X):=\{(z,w)\in X\times \C^m: H(z,w)<1\}
$$
and call it a \emph{Hartogs type} domain \cite{TM}.  Hartogs domains correspond to the special case
$m=1$ and $H(z,w)= |w|e^{\varphi(z)}$.

Motivated by studying hyperbolicity and tautness modulo an analytic subset of complex spaces, the main goal of this article is to give necessary and sufficient conditions on hyperbolicity or tautness modulo a ``vertical" analytic subset of the Hartogs domains $\Omega_{\varphi}(X).$ The results are given in Section 
\ref{results}, but first we recall some basic notions.

\begin{definition} (see \cite[p. 68]{K}) \ Let $X$ be a complex space and $S$ be an analytic subset of $X.$  We say that $X$ is \emph{hyperbolic modulo $S$} if for every pair of distinct  points $p,q$ of $X$ we have $d_X(p,q)>0$ unless both are contained in $S,$ where $d_X$ is the Kobayashi pseudodistance of $X.$ 
\end{definition}
If $S=\emptyset,$ then $X$ is said to be hyperbolic.

\begin{definition} (see \cite [p.240]{K}) \ Let $X$ be a complex space and $S$ be an analytic subset in $X$. We say that $X$ is \emph{taut modulo $S$} if it is 
normal modulo $S$, i.e., for every sequence $\{f_n\}$ in $Hol(\mathbb D,X)$ one of the following holds:
\begin{enumerate}
\item[i.] There exists a subsequence of $\{f_n\}$ which converges uniformly to $f\in Hol(\mathbb D,X)$ in $Hol(\mathbb D,X)$;
\item[ii.] The sequence $\{f_n\}$ is compactly divergent modulo $S$ in $Hol(\mathbb D,X)$, i.e., for each compact set $K\subset\mathbb D$ and each compact set $L\subset X\setminus S$, there exists an integer $N$ such that $f_n(K)\cap L=\emptyset$ for all $n\geq N.$
\end{enumerate}
\end{definition}
If $S=\emptyset,$ then $X$ is said to be taut.  It is immediate from the definition
that if $S \subset S' \subset X$ and $X$ is taut modulo $S$, then it is taut modulo $S'$,
so in particular if $X$ is taut, it is taut modulo $S$ for any analytic subset $S$. 

Of course, the converse does not hold. 

\begin{example} 
\label{ex1}
Let $X=\{(z,w)\in \mathbb C^2| |z|<1, |zw|<1\}$ and $S:=\{0\}\times \mathbb C.$  Then 
\begin{enumerate}
\item[(i)] $X$ is not hyperbolic, but $X$ is hyperbolic modulo $S.$
\item[(ii)] $X$ is not taut, but $X$ is taut modulo $S.$
\item[(iii)] $X\setminus S$ is taut (thus hyperbolic).
\end{enumerate}
\end{example}

We could have $X\setminus S$ being taut without $X$ being taut modulo $S$: for instance,
$\C\setminus \{0,1\}$ is taut, but $\C$ is not taut modulo $\{0,1\}$.
 
On the other hand, there are examples of domains taut modulo $S$ such that $X \setminus S$
is not taut: just take $X$ a taut domain and $S$ such that the codimension of $ S$ is at least $ 2$. Then $X\setminus S$ is not pseudoconvex, therefore not taut. 

\begin{proof*}{\it Proof of Example \ref{ex1}.}
Since the complex line $S$ is contained in $X$, it cannot be hyperbolic, thus 
isn't taut either.  On the other hand, $X\setminus S$ is biholomorphic 
to $(\D\setminus\{0\})\times \D$ under the map $(z,w) \mapsto (z,zw)$,
and the latter is clearly taut, which proves (iii).  

Now suppose $(z_0,w_0) \neq (z_1,w_1) \in X$, with at least
one of them not in $S$. If $z_0\neq z_1$, then 
$d_X \left( (z_0,w_0), (z_1,w_1) \right) \ge d_{\D} (z_0,z_1) >0$; if $z_0= z_1$,
then $z_0\neq 0$.  Given any finite set of points of $X$, $(\zeta_k, \eta_k)$
which connect $(z_0,w_0)$ to $(z_1,w_1)$ via consecutive analytic disks, either
there is some $k$ such that $|\zeta_k| \le |z_0|/2$ and then the corresponding sum
will contribute at least $d_{\D} (z_0,z_0/2)$, or there is not, and then
all points are in $X \cap \{ |z| > |z_0|/2\} \subset \{  |z_0|/2 < |z| <1, |w| < 2/|z_0|\}:=P$ and  the sum
will be bounded below by $d_P \left( (z_0,w_0), (z_1,w_1) \right) \ge \frac2{|z_0|} d_{\D}(w_0,w_1)>0$. So $X$ is hyperbolic modulo $S.$

\noindent
{\it Proof of {(ii)}.} Assume that $\{f_n\}\subset Hol(\mathbb D, X)$ is not compactly divergent modulo $S$. Write $f_n=(g_n,h_n)$ with $g_n, h_n\in Hol(\mathbb D, \mathbb C)$ satisfying $|g_n(z)|<1$, $|g_n(z)h_n(z)|<1$ for all $z\in \mathbb D$ and for $n=1,2,\cdots$ Since $\{f_n\}$ is not compactly divergent modulo $S$, $\{g_n\}$ is not compactly divergent modulo $\{0\}$. Therefore, by Montel's theorem we may assume, without loss of generality, that $\{g_n\}$ converges uniformly on every compact subset of $\mathbb D$ to a holomorphic function $g\in Hol(\mathbb D,\mathbb D)$, not identically zero. Since $\{g_nh_n\}$ is not compactly divergent on $\mathbb D$ and since $\mathbb D$ is taut, taking a subsequence we may assume that $\{g_nh_n\}$ also converges uniformly on every compact subset of $\mathbb D$ to a holomorphic function $\gamma\in Hol(\mathbb D,\mathbb C)$. Hence $\{h_n\}$ converges uniformly on every compact subset of $\mathbb D$ to a meromorphic function $h:= \gamma/g$ on $\mathbb D$. Moreover, by Hurwitz's theorem $h$ is actually holomorphic on $\mathbb D$ and thus $\{f_n\}$ converges uniformly on every compact subset of $\mathbb D$ to a holomorphic map $f:=(g,h)\in Hol(\mathbb D,\overline{X})$.

We now prove that $f\in Hol(\mathbb D,X) $. Since $g\in Hol(\mathbb D, \mathbb D)$, it suffices to show that $|g(z)h(z)|<1$ for all $z\in\mathbb D$. Indeed, suppose not. Then there is $z_0\in\mathbb D$ such that $|g(z_0)h(z_0)|=1$. By the maximum principle, $gh$ is a constant function. Therefore, $|g(z)h(z)|=1$ for every $z\in\mathbb D$. This is not possible because $\{g_nh_n\}\subset Hol(\mathbb D,\mathbb D)$ is not compactly divergent. Thus, the proof is complete.
\end{proof*}

\section{Main Results}
\label{results}
We denote $\widetilde{S}:=S\times\C^m$. Recall that $H: X\times \C^m \to [-\infty,\infty)$
is an upper semicontinuous function such that $H(z,w) \ge 0, H(z,\lambda w)=|\lambda|H(z,w), \lambda \in \C, z \in X, w \in \C^m$.

\begin{theorem}
\label{hyperbolicity} 
Let $X$ be a complex space and $S$ be an analytic subset in $X$. 
Then $\Omega_H(X)$ is hyperbolic modulo $\widetilde{S}$ if and only if $X$ is hyperbolic modulo $S$ and the function $H$
 satisfies the following condition:
\begin{multline}
\label{hypcond}
\mbox{ If }\{z_k\}_{k \ge 1} \subset X\setminus S\mbox{ with }\lim \limits_{k \to \infty}z_k = z_0\in X\setminus S\\
\mbox{ and  }
\{w_k\}_{k \ge 1}\subset \C^m\mbox{ with }\lim\limits_{k \to \infty}w_k = w_0 \not= 0
\mbox{, then }\liminf\limits_{k \to \infty}H(z_k,w_k) \not= 0 .
\end{multline}
\end{theorem}

The proof is given in Section \ref{hyper}.

\begin{corollary}
Let $X$ be a complex space, $S$ an analytic subset in $X$, $\varphi: X\to [-\infty,\infty)$ an upper-semicontinuous function on $X$.
Then the Hartogs domain
$\Omega_\varphi(X)$ is hyperbolic modulo $\widetilde{S}$ if and only if $X$ is hyperbolic modulo $S$ and $\varphi$ is locally bounded (below) on $X\setminus S$. 
\end{corollary}

The situation for tautness is a bit more complicated, at least in the case of complex spaces. 

\begin{theorem}\label{Thm2} 
Let $X$ be a complex space and $S$ be an analytic subset in $X$. 
Then 
\begin{enumerate}
\item[i.] If  
$\Omega_H(X)$ is taut modulo $\widetilde{S}$, then $X$ is taut modulo $S$ and $\log H$ is continuous plurisubharmonic on $X\setminus S \times \C^m$.
\item[ii.] If furthermore $X$ is a complex manifold and $S$ a (proper) analytic subset,
then in addition $\log H$ is  plurisubharmonic on $X \times \C^m $.
\item[iii.] Conversely, if $X$ is taut modulo $S, H$ is continuous on $X\setminus S \times \C^m$ and  $\log H$ is plurisubharmonic on $X \times \C^m,$ then $\Omega_H(X)$ is taut modulo $\widetilde{S}.$ 
\end{enumerate}
\end{theorem}

The proof is given in Section \ref{taut}.

As above, an immediate corollary is obtained for Hartogs domains by observing
that $H(z,w):= |w|e^{\varphi(z)}$ is continuous if and only if $\varphi$ is,
and
$\log H$ is plurisubharmonic if and only if $\varphi$ is. 

Implication (ii) cannot hold as stated for the general case of a complex space, as the following shows.

\begin{example}
\label{notsh}
Let $X:=\{ (z_1,z_2)\in\C^2: z_1z_2=0\}$, $S:= \{ (z_1,z_2)\in\C^2: z_2=0\}$,
$\varphi(z_1,z_2):= \log |z_2|$. Then $\Omega_\varphi(X)$ is taut modulo $S$,
but $\varphi \notin PSH(X)$.
\end{example}

\begin{proof}
Since $\varphi$ is identically $-\infty$ on an open set, it can't be plurisubharmonic.
If a sequence $(F_n)\subset Hol(\D,\Omega_\varphi(X))$ is not compactly divergent modulo $\widetilde S$, 
it is easy 
to see that $F_n(\D) \subset \{(z,w) \in \Omega_\varphi(X) : z_1=0\}$. But then if $F_n=(0,g_n,h_n)$, its convergence is equivalent to that of 
$(g_n,h_n)\subset Hol(\D,\Omega_{\varphi_0}(\C))$, where $\varphi_0(z_2)= \log |z_2|$.
We know that $(g_n,h_n)$ is not compactly divergent modulo $\{0\}$, and Theorem
\ref{Thm2} shows that $\Omega_{\varphi_0}(\C)$ is taut modulo $\{0\}$, so a subsequence must
converge on compacta.
\end{proof}

Note that in this case $X$ is not irreducible and $S$ is a whole component. We don't 
know what happens if we rule out this degenerate situation.

\section{Hyperbolicity}
\label{hyper}
%{Proof of Theorem \ref{hyperbolicity}}
 
Recall that the Lempert function is defined as
$\ell_{\Omega}(a,b)=\inf\{p(0,\lambda): \exists \varphi\in Hol(\D,\Omega), \varphi(0)=a, \varphi(\lambda)=b\}$ and that $\ell_{\Omega}(a,b) \le k_{\Omega}(a,b)$.

Using the same argument as in the proof of Remark 3.1.7 and Proposition 3.1.10 in \cite{JP},
 it is easy to see that

\begin{lemma}
\label{poincare}
Let $\Omega=\Omega_H(X).$ Then $\ell_{\Omega}((z,0),(z,w))\le p(0,H(z,w))$ for any $(z,w)\in \Omega,$
 where $p$ is the Poincar\'e distance. 
Equality holds if $H \in PSH(X\times \C^m).$ 
\end{lemma}
 
\noindent
\begin{proof*}{\bf Proof of Theorem \ref{hyperbolicity}}

\noindent
$(\Longrightarrow)\ $ Suppose $\Omega_H(X)$ is hyperbolic. Since $X$ is isomorphic to
 a closed complex subspace of $\Omega_H(X)$, we deduce that $X$ is hyperbolic. 
Next, we will show that $H$ verifies the property \eqref{hypcond}. Otherwise, there would exist
 $ \{z_k\}_{k \ge 1} \subset X\setminus S$ with $ \lim \limits_{k \to \infty}z_k = z_0\in X \setminus  S,  \{w_k\}_{k \ge 1}
\subset \C^m$ with $\lim\limits_{k \to \infty}w_k = w_0 \not= 0$ such that $\lim\limits_{k \to \infty}H(z_k,w_k) = 0.$
 Without loss of generality, we may assume that $(z_k,w_k)\in\Omega_H(X).$ 
Then by Lemma \ref{poincare}, we have 
$$
0 \le k_{\Omega}((z_k,0),(z_k,w_k)) \le p(0,H(z_k,w_k)), \forall k \ge 1.
$$
By letting $k$ go to $\infty$, we find that $k_{\Omega}((z_0,0),(z_0,w_0))=0.$
 This contradicts  the hyperbolicity modulo $\tilde S$ of $\Omega_H(X)$.

\noindent 
$(\Longleftarrow)$ To prove the converse, we consider the projection $\pi: \Omega_H(X) \to X$ 
given by $\pi(z,w)=z.$ Let $U$ be a compact neighbourhood of $z_0$ in $X \setminus S.$
 Then $ \underset{z\in U}{\cup}\Omega_H(z)$ is a bounded set in $\C^m,$ where $ \Omega_H(z):=\{w\in\C^m:H(z,w)<1\}.$ 

In fact, suppose that this property does not hold. Then $\exists \{z_k\}_{k\ge 1} \subset U, \{w_k\}_{k\ge 1}\subset \C^m$
 such that $\lim \limits_{k \to \infty}\|w_k\|=\infty$ and $H(z_k,w_k)<1.$ Put $w_k:=r_ku_k$ 
with $\|u_k\|=1, \forall k\ge1$; then $|r_k|\to \infty$ as $k \to \infty.$ Passing to a subsequence,
 we may assume that $z_k \to z_0$ and $u_k \to u_0 \not= 0$ as $k \to \infty.$ 
Since $H(z_k,w_k)=|r_k|H(z_k,u_k)<1,$  we have $\limsup\limits_{k\to\infty}H(z_k,u_k)=0.$ 
This contradicts the property \eqref{hypcond}. So, there exists $R>0$ 
such that  
$\pi^{-1}(U)\subset U\times \underset{z\in U}{\cup}\Omega_H(z)\subset U\times B(0,R).$ 
Therefore, $\pi^{-1}(U)$ is hyperbolic too. By Eastwood's theorem  \cite{Ea} we conclude the proof.
\end{proof*}

\section{Tautness}
\label{taut}

\begin{proof*}{\bf Proof of Theorem \ref{Thm2}}

\noindent {\it Proof of (i).} 

 Since $X$ is isomorphic to a closed complex subspace of $\Omega_H(X)$, we deduce that $X$ is taut modulo $S.$ 
We now show that $H$ is continuous on $X\setminus S \times \C^m.$ Otherwise, there would exist
 $r >0, \{(z_k,w_k)\}_{k \ge 1} \subset X\setminus S \times \C^m$ such that
$$
\{(z_k,w_k)\} \to (z_0,w_0), z_0\in X\setminus S  \text { and }  H(z_k,w_k)<r<H(z_0,w_0), \forall k \ge 1.
$$

For each $k \ge 1,$ we define the holomorphic mapping $f_k: \D \to \Omega_H(X)$ given by $f_k(\lambda)=(z_k,\displaystyle{\frac{\lambda w_k} r}).$ 
It is clear that $f_k(0)=(z_k,0) \to (z_0,0) \in \Omega_H(X)\setminus \tilde S.$ Since $\Omega_H(X)$ is taut modulo $\tilde S$, 
by passing to a subsequence if necessary, we may assume that $f_k$ converges locally uniformly on $\D$ 
to a holomorphic mapping $f \in Hol(\D,\Omega_H(X)).$ It is easy to see that 
$f(\lambda)=(z_0,\displaystyle{\frac{\lambda w_0} r}).$ Hence
$$
\displaystyle{\frac{|\lambda|} r}H(z_0,w_0)=H(z_0,\displaystyle{\frac{\lambda w_0} r})<1, \forall \lambda \in \D.
$$
This implies that $H(z_0,w_0)<\displaystyle{\frac{r}{|\lambda|}}, \forall \lambda \in \D,$
and hence $H(z_0,w_0) \le r.$ This is a contradiction.

It remains to show that $\log H$ is plurisubharmonic. 

According to the theorem of Fornaess and Narasimhan \cite{FN}, it suffices to show that 
$u(z):=\log H \circ g(z)=\log H(g_1(z),g_2(z))$ is subharmonic for every $g =(g_1,g_2)\in Hol(\D,X\setminus S\times \C^m)\cap C(\overline{\D},X\setminus S\times \C^m).$ Suppose the contrary. Then $\exists z_0 \in \D, r > 0$ such that $\overline D(z_0,r) \subset \D$ and a harmonic function $h$ such that
$h(z)\ge u(z)$ for any $z=z_0+re^{i\theta}, \forall \theta \in \R,$ but $u(z_0)>h(z_0).$ 
Let $\tilde h$ denote a harmonic conjugate to $h$. 

We have $u(z)-h(z)=\log H\Big(g_1(z),e^{-h(z)-i\tilde{h}(z)}g_2(z)\Big) \le 0, \forall z=z_0+re^{i\theta}$
and $u(z_0)-h(z_0)=\varepsilon_0>0.$

For any $n \ge 1,$ we set $\varphi_n(\lambda):=\Big(g_1(z),e^{-h(z)-i\tilde{h}(z)-\varepsilon_0 -\frac 1 n}g_2(z)\Big),$ 
where $z\in \overline{D}(z_0,r).$ Then
 $\varphi_n \in Hol(D(z_0,r),\Omega_H(X)) \cap C(\overline{D}(z_0,r),\Omega_H(X))$, 
 $\cup_{n=1}^{\infty}\varphi_n(\partial D(z_0,r))\Subset \Omega_H(X),$ 
and $\varphi_n(0)$ tend to a boundary point. 
This contradicts  the tautness of $\Omega_H(X).$

\noindent
{\it Proof of (ii).}
By the removable singularity theorem for plurisubharmonic functions
 \cite[Theorem 2.9.22]{Kl} since $\log H$ is locally bounded above in $X\times \C^m$ and plurishubharmonic
in $X\setminus S \times \C^m$, it can be extended across $\tilde S$
to a function $\log \hat H \in PSH(X\times \C^m)$  given, for $(z_0,w_0)\in \tilde S$, by
$$
\log \hat H (z_0,w_0) := \limsup_{(z,w)\to (z_0,w_0), (z,w) \notin \tilde S} H(z,w).
$$
We claim that $\hat H = H$.  Since $H$ is upper semi continuous
on $X\times \C^m$, this conclusion can only fail if there exists $(z_0,w_0) \in \tilde S$ such that 
$H(z_0,w_0) > \limsup_{(z,w)\to (z_0,w_0), (z,w) \notin \tilde S} H(z,w)$.  Since $X$ is a manifold,
we can go to a coordinate patch and find an analytic disk $f$ such that 
$f(0)=(z_0,w_0)$, $f(\D) \not \subset \tilde S$.  Then $f^{-1}(\tilde S)$ must be a discrete subset 
of $\D$, and reducing the disk we may assume that $f^{-1}(\tilde S)= \{0\}$ and
$\sup_{0<|\zeta|} H(f(\zeta))= H(z_0,w_0)-\delta$, $\delta>0$.  The 
proof then proceeds essentially as above. We have a contradiction.

\noindent
%($\Leftarrow$) 
{\it Proof of (iii).} Assume that $X$ is taut modulo $S$ and $\log H$ is continuous on $X\setminus S\times \C^m$, and plurisubharmonic on $X\times\C^m$. We now show that $\Omega_H(X)$ is taut modulo $\widetilde{S}$.

Consider the projection $\pi:\Omega_H(X)\rightarrow X$ defined by $\pi(x,z)=x$. We now prove, for each $x\in X\setminus S$, that there exists an open neighbourhood $U$ of $x$ in $X\setminus S$ such that $\pi^{-1}(U)$ is taut. Indeed, choose a hyperconvex neighbourhood $U$ of $x$ in $X\setminus S$. It is easy to see that $\pi^{-1}(U)=\{(u,z)\in U\times\C^m: H(z,w)<1\}=\Omega_H(U).$ Suppose that $\rho$ is a negative plurisubharmonic exhaustion function of $U$. Then $(u,z)\mapsto\max(\rho(u),\log H(u,z))$ is also a negative plurisubharmonic exhaustion function of $\Omega_H(U)$. Thus, $\Omega_H(U)$ is hyperconvex. By a Theorem of Sibony \cite{Si} 
and \cite{TD}, $\Omega_H(U)$ is taut. Thus,  $\pi^{-1}(U)$ is taut. 

Assume that a sequence  $\{\widetilde{f}_n\}\subset Hol(\mathbb D,\Omega_H(X))$ is not compactly divergent modulo $\widetilde{S}$
in $Hol(\mathbb D,\Omega_H(X)).$ From now on, denote
$\widetilde{X}:=\Omega_H(X)\setminus \widetilde{S}$.

Without loss of generality, 
we may assume that there exist a compact set $K\subset \mathbb D$ and a compact set 
$L\subset \widetilde{X}$ 
such that $\widetilde{f}_n(K)\cap L\not=\emptyset$ for $n\geq 1$. 
For each $n\geq 1$, there exists $z_n\in K\subset \mathbb D$ such that $\widetilde{f}_n(z_n)\in L$. Since $K$ and $L$ are compact sets, by taking subsequences if necessary, 
we may assume that $\{z_n\}\subset K\subset \mathbb D$ such that 
$z_n\rightarrow z_\infty\in K\subset \mathbb D$ and 
$\widetilde{f}_n(z_n)\rightarrow \widetilde{p}\in L\subset\widetilde{X}$.
 It is easy to see that $\{f_n:=\pi\circ \widetilde{f}_n\}\subset Hol(\mathbb D,X)$ is not compactly divergent modulo $S$ in $Hol(\mathbb D,X)$. Since $X$ is taut modulo $S$, we may assume that $\{f_n\}$ converges uniformly to a mapping $F\in Hol(\mathbb D,X)$. Obviously, $\pi(\widetilde{f}_n(z_n))\rightarrow\pi(\widetilde{p})$ and $\pi(\widetilde{f}_n(z_n))=\pi\circ \widetilde{f}_n(z_n)=f_n(z_n)\rightarrow F(z_\infty)$ as $n\to\infty$. Therefore, we can let $p=\pi(\widetilde{p})=F(z_\infty)$.

Since $\widetilde{p}\in\widetilde{X}$, $p=\pi(\widetilde{p})\not\in S$. 
Then there exists an open neighbourhood $U$ of $p$ in $X\setminus S$ 
such that $\pi^{-1}(U)$ is taut. 
Taking an open neighbourhood $V\Subset F^{-1}(U)$ of $z_\infty$ in $\mathbb D\setminus F^{-1}(S)$ 
and since the sequence $\{f_n\}$ converges uniformly to a mapping $F$, 
we may assume that  $f_n(V)\subset U$. 
This implies that $\widetilde{f}_n(V)\subset\pi^{-1}(U)$ for every $n\geq 1$.

Consider the compact subsets $K=\{z_n, n\in \Z_+\}\cup \{z_\infty\}\subset \mathbb D$ and
 $L=\{\widetilde{f}_n(z_n)\}\cup \{\widetilde{p}\}\subset \widetilde X$. 
 Then $\widetilde{f}_n(K)\cap L\not=\emptyset$ for all $n,$ and hence the sequence $\{\widetilde{f}_n\bigr|_V\}$ is not compactly divergent modulo $\widetilde{S}$ in $Hol(\mathbb D,\Omega_H(X))$. 
 Since $\pi^{-1}(U)$ is taut and $\widetilde{f}_n(V)\subset\pi^{-1}(U)$, 
 it implies that $\{\widetilde{f}_n\bigr|_V\}$ converges uniformly to a mapping $\widetilde{F}$ in $Hol(V,\pi^{-1}(U))$.

Consider the family $\Gamma$ of all pairs $(W,\Phi)$, where $W$ is an open set in $\mathbb D\setminus F^{-1}(S)$ and $\Phi\in Hol(W,\widetilde{X})$ 
such that there exists a subsequence $\{\widetilde{f}_{n_k}\bigr|_W\}$ of $\{\widetilde{f}_{n}\bigr|_W\}$ 
which converges uniformly to mapping $\Phi$ in $Hol(W,\widetilde{X}).$

According to the proof above, we have $ \Gamma\not=\emptyset.$ We now consider the following order relation  in the family $\Gamma$: $(W_1,\Phi_1)\leq (W_2,\Phi_2)$ if

i) $W_1\subset W_2$ and

ii)  for any subsequence $\{\widetilde{f}_{n_k}\bigr|_{W_1}\}$ of $\{\widetilde{f}_n\bigr|_{W_1}\}$ that converges uniformly to mapping $\Phi_1$ in $Hol(W_1,\widetilde{X})$, there exists a subsequence $\{\widetilde{f}_{n_{k_l}}\}$ of $\{\widetilde{f}_{n_k}\}$ such that the sequence $\{\widetilde{f}_{n_{k_l}}\bigr|_{W_2}\}$ converges uniformly to mapping $\Phi_2$ in $Hol(W_2,\widetilde{X})$.

Assume that $\{(W_\alpha,\Phi_\alpha)\}_{\alpha\in\Lambda}$ is a well-ordered subset of $\Gamma$.
Put  $W_0=\bigcup\limits_{\alpha\in\Lambda}^{}W_\alpha$ and define a mapping $\Phi_0\in Hol(W_0,\Omega_H(X)\setminus\widetilde{S})$ given by $\Phi_0\bigr|_{W_\alpha}=\Phi_\alpha$ for each $\alpha \in \Lambda$. 

Take a sequence $\{(W_i,\Phi_i)\}_{i=1}^{\infty}\subset \{(W_\alpha,\Phi_\alpha)\}_{\alpha\in\Lambda}$ such that
$$(W_1,\Phi_1)\leq (W_2,\Phi_2)\leq\cdots \text{ and } W_0=\bigcup\limits_{i=1}^{\infty}W_i.$$
By the definition of $\Gamma$, there exists a subsequence $\{\widetilde{f}_{n}^1\bigr|_{W_1}\}$ of $\{\widetilde{f}_n\bigr|_{W_1}\}$ such
that  $\{\widetilde{f}_{n}^1\bigr|_{W_1}\}$  converges uniformly to the mapping $\Phi_1$ in 
$Hol(W_1,\widetilde{X}).$ By the definition of the order relation on $\Gamma$, there exists a subsequence $\{\widetilde{f}_n^2\}$ of $\{\widetilde{f}_n^1\}$ which converges uniformly to $\Phi_2$ in $Hol(W_2,\widetilde{X}).$

By continuing this process, we get the sequence $\{\widetilde{f}_n^k\}$ such that $\{\widetilde{f}_n^k\}\subset\{\widetilde{f}_n^{k-1}\}$ for all $k\geq 2$ and $\{\widetilde{f}_n^k\bigr|_{W_k}\}$ converges uniformly to $\Phi_k$ in $Hol(W_k,\widetilde{X})$. Thus, a diagonal sequence $\{\widetilde{f}_k^k\}$ converges uniformly to $\Phi_0$ in $Hol(W_0,\widetilde{X})$. Hence $(W_0,\Phi_0)\in\Gamma$ and the subset $\{(W_\alpha,\Phi_\alpha)\}_{\alpha\in\Lambda}$ of $\Gamma$ has a supremum.

By the Zorn lemma, there exists a maximal element $(W,\Phi)$ of the family $\Gamma$. Assume that $\{\widetilde{f}_{n_k}\bigr|_W\}$ is a subsequence of $\{\widetilde{f}_{n}\bigr|_W\}$ such that $\{\widetilde{f}_{n_k}\bigr|_W\}$ converges uniformly to $\Phi$ in $Hol(W,\widetilde{X})$. 

We now show that $W=\mathbb D\setminus F^{-1}(S)$. Suppose that there exists $z_0\in \overline W\cap( \mathbb D\setminus F^{-1}(S))$. Take an open neighbourhood $U$ of $F(z_0)$ in $X\setminus S$ such that $\pi^{-1}(U)$ is taut. Since $\{f_n\}$ converges uniformly to a mapping $F$ in $Hol(\mathbb D,X)$, there exists an open neighbourhood $W_0$ of $z_0$ in $\mathbb D\setminus F^{-1}(S)$ such that $\pi\circ\widetilde{f}_n(W_0)\subset U$. Hence $\widetilde{f}_n(W_0)\subset \pi^{-1}(U)$  for all $n\geq 1$.

Fix $z_1\in W_0\cap W$. Then the sequence $\{\widetilde{f}_{n_k}(z_1)\}$ is convergent. Since $Hol(W_0,\pi^{-1}(U))$ is a normal family, $\{\widetilde{f}_{n_k}\bigr|_{W_0}\}$ converges uniformly to $\Phi_0$ in $Hol(W_0,\pi^{-1}(U))$. Thus $(W_0,\Phi_0)\in\Gamma$. It implies that $W_0\subset W$. Hence $W=\mathbb D\setminus F^{-1}(S)$.

 Since $F^{-1}(S)$ is an analytic subset in the open unit disc $\mathbb D$, $F^{-1}(S)$ is a discrete set and hence, $F^{-1}(S)$ does not have 
any accumulation point in $\mathbb D$. Therefore, for each $z\in F^{-1}(S)$, there exists a number $0<r<1$ such that $B(z,r)\cap F^{-1}(S)=\{z\}$. Without loss of generality, we may assume that $\mathbb D\setminus F^{-1}(S)=\mathbb D\setminus \{0\}:=\mathbb D^*$  and the sequence $\{\widetilde{f}_n\bigr|_{\mathbb D^*}\}$  converges uniformly to the map $\Phi\in Hol(\mathbb D^*, \widetilde X)$.

We  now rewrite functions $\tilde f_n$ and $\Phi$ as follows: $\tilde f_n(z)=(f_n(z),g_n(z))$ for each $z\in \mathbb D$ and $\Phi(z)=(F(z),G(z))$ for each $z\in \mathbb D^*$. Since the sequence $\{g_n\}$ converges uniformly on compact subsets of $\mathbb D^*$ to $G\in Hol(\mathbb D^*,\mathbb C^m)$, the maximum principle implies that it is a uniformly
Cauchy sequence on compact subsets of $\D$, and therefore converges uniformly  on compact subsets of $\mathbb D$ to $\tilde G\in Hol(\mathbb D,\mathbb C^m)$. Denote by $\tilde \Phi=(F,\tilde G): \mathbb D \to X\times \mathbb C^m$. We now prove that $\tilde \Phi\in Hol(\mathbb D,\Omega_H(X))$. Since $\tilde \Phi(\mathbb D^*)=\Phi(\mathbb D^*)\subset \Omega_H(X)$, it suffices to show that $\tilde \Phi(0)\in \Omega_H(X)$. 

Indeed,  since $H$ is plurisubharmonic on $X\times \C^m$, $H \circ \tilde \Phi$
is subharmonic on $\D$, and since $\Phi(\mathbb D^*)\subset \Omega_H(X)$, it is 
negative on $\D^*$. So it is negative on $\D$, i.e. $\tilde \Phi (\mathbb D)\subset \Omega_H(X)$. Hence $\Omega_H(X)$ is taut modulo $\widetilde{S}$.
\end{proof*}

\bibliographystyle{amsplain}

\begin{thebibliography}{99}
\vspace{20pt}
\bibitem{A} Y. Adachi, {\it On a Kobayashi hyperbolic manifold N modulo a closed subset   and its applications}, 
Kodai Math. J. 30 (2007), 131--139.
\bibitem{DT} Nguyen Quang Dieu and Do Duc Thai, {\it Complete hyperbolicity of Hartogs domain}, Manuscripta Math. 112 (2003), 
171--181.
\bibitem{Ea} A. Eastwood, {\it \`A propos des vari\'et\'es hyperboliques compl\`etes}, C. R. Acad. Sci. Paris {\bf 280} (1975), 1071--1075.
\bibitem{FN} J. Fornaess and R. Narasimhan, {\it The Levi problem on complex spaces with singularities}, Math. Ann. 248 (1980), 47--72.
\bibitem{JP} M. Jarnicki and P. Pflug, {\it Invariant Distances and Metrics in Complex Analysis}, Walter de Gruyter, Berlin-New York (1993).
\bibitem{Kl} M. Klimek, {\it Pluripotential Theory}, Oxford University Press, New York, 1991.
\bibitem{K} S. Kobayashi, {\it Hyperbolic complex spaces}, Springer-Verlag, Berlin, 1998.
\bibitem{R} R. Remmert, {\it Classical topics in complex function theory}, Springer-Verlag, Berlin, 1998.
\bibitem{Si} N. Sibony, {\it A class of hyperbolic manifolds}, Recent Developments in Several Complex Variables, Princeton Univ. Press 100 (1981), 357--372.
\bibitem{TT} Do Duc Thai and P. J. Thomas, {\it $D^*$-extension property without hyperbolicity}, Indiana Univ. Math. J. 47 (1998), 1125--1130.
\bibitem{TD} Do Duc Thai and Pham Viet Duc, {\it On the complete hyperbolicity and the tautness of the Hartogs domains}, 
Intern. J. Math. 11 (2000), 103--111.
\bibitem{TM} Nguyen Van Trao and Tran Hue Minh, {\it  Remarks on the Kobayashi hyperbolicity of complex spaces}, Acta Math. Vietnam. 
34 (2009),  375--387.

\vskip0.5cm
\end{thebibliography}

\end{document}